# Some Optimization Problems with Calculus.

## Arkady Kitover and Mehmet Orhon

**Abstract.** Starting from the well-known and elementary problem of inscribing the rectangle of the greatest area in an ellipse, we look at gradually more and more complicated variants of this problem. Our goal is to demonstrate to an average but motivated student of Calculus how to start from an inconspicuous textbook problem and to arrive at considerably more interesting and complicated problems, some of which can be subjects of independent research.

**Introduction.** The following problem is a standard staple of Calculus textbooks (see e.g. [SCW, p. 344, Problem 30]). What are the dimensions of the rectangle of the greatest area inscribed in the ellipse $\frac{x^2}{a^2} + \frac{y^2}{b^2} = 1$? A little bit more interesting (but also simple) problem of finding the rectangle of the greatest perimeter inscribed in the ellipse is employed in Calculus textbooks considerably less often. We add to these two problems the third one: to find the rectangle inscribed in the ellipse with the greatest ratio $\frac{S}{P^2}$, where $S$ is the area of the rectangle and $P$ is its perimeter.

After looking at the solutions to these problems in the case of an ellipse we discuss possible generalizations in three incremental steps.

I. We consider rectangles inscribed in the curve $\frac{|x|^\alpha}{A^\alpha} + \frac{|y|^\alpha}{B^\alpha} = 1$, where $A, B,$ and $\alpha$ are positive real numbers.

II. We slightly change the problem by considering the rectangles inscribed in the curve $\frac{|x|^\alpha}{A^\alpha} + \frac{|y|^\beta}{B^\beta} = 1, \alpha > 0, \beta > 0$. As the reader will see, this small change makes some of our three problems more difficult to solve.

III. We discuss the case of differentiable (or piecewise differentiable) closed curves symmetric about the $x$ and $y$ - axes, state some results, and pose some open problems. [1]

In the second part of the paper, we consider similar problems related to rectangular parallelepipeds inscribed in surfaces symmetric about the coordinate planes. As can be expected, some of these problems are considerably more interesting and difficult than the corresponding two-dimensional problems.

## Part 1. Optimization problems in two dimensions.

### Part 1A. Ellipse.

(a) Find the dimensions of the rectangle of the greatest area inscribed in the ellipse $\frac{x^2}{a^2} + \frac{y^2}{b^2} = 1$. Notice that the area of an inscribed rectangle is $4xy$, where $(x, y)$ is the vertex of the rectangle in the first quadrant. This expression attains its maximum together with the product $\frac{x^2}{a^2} \cdot \frac{y^2}{b^2}$. Using the well-known fact that the product of two nonnegative numbers $u$ and $v$ with a fixed sum $s$ takes its greatest value if and only if $u = v = \frac{s}{2}$ we see that $\frac{x^2}{a^2} = \frac{y^2}{b^2} = \frac{1}{2}$. Hence, the dimensions of the rectangle of the greatest area inscribed in the ellipse are $2x = a\sqrt{2}, 2y = b\sqrt{2}$ and its area is $2ab$.

**Remark 1.** The rectangle with the area $2ab$ inscribed in the ellipse $\frac{x^2}{a^2} + \frac{y^2}{b^2} = 1$ is unique, but there are infinitely many (continuum) quadrilaterals of the same area inscribed in this ellipse. All of them are parallelograms and one of them is the rhombus with the vertices $(a,0),(-a,0),(0,b),(0,-b)$. We refer the reader to [Jo] for an excellent presentation of proof of these facts. [2]

---

[1] We would like to emphasize that optimization is not our area of expertise, and when we state that a problem is open, we mean by it that we were not able to find its solution in the literature.

(b) Find the dimensions of the rectangle of the greatest perimeter inscribed in the ellipse $\frac{x^2}{a^2}+\frac{y^2}{b^2}=1$. We have to maximize the sum $x+y$ subject to $\frac{x^2}{a^2}+\frac{y^2}{b^2}=1, x>0, y>0$. The method of Lagrange multipliers provides the equations,

$$1=\frac{2\lambda x}{a^2}=\frac{2\lambda y}{b^2}.$$

Thus, $y=\frac{b^2}{a^2}x$, and combining this equation with the equation of the ellipse we easily conclude that,

$$x=\frac{a^2}{\sqrt{a^2+b^2}}, y=\frac{b^2}{\sqrt{a^2+b^2}}.$$

The maximum perimeter is $P_{max}=4\sqrt{a^2+b^2}$. Notice that the rhombus inscribed in the ellipse has the same perimeter. Figure 1 below shows the rectangle of the greatest perimeter and the rhombus inscribed in the ellipse $\frac{x^2}{16}+\frac{y^2}{9}=1$.

---

[2] To avoid possible ambiguity let us agree that by a rectangle inscribed in a curve symmetric about the coordinate axes we always understand a rectangle with the vertices on the curve and the axes of symmetry parallel to the coordinate axes. Similar agreement relates to rectangular parallelepipeds inscribed in surfaces symmetric about coordinate planes.

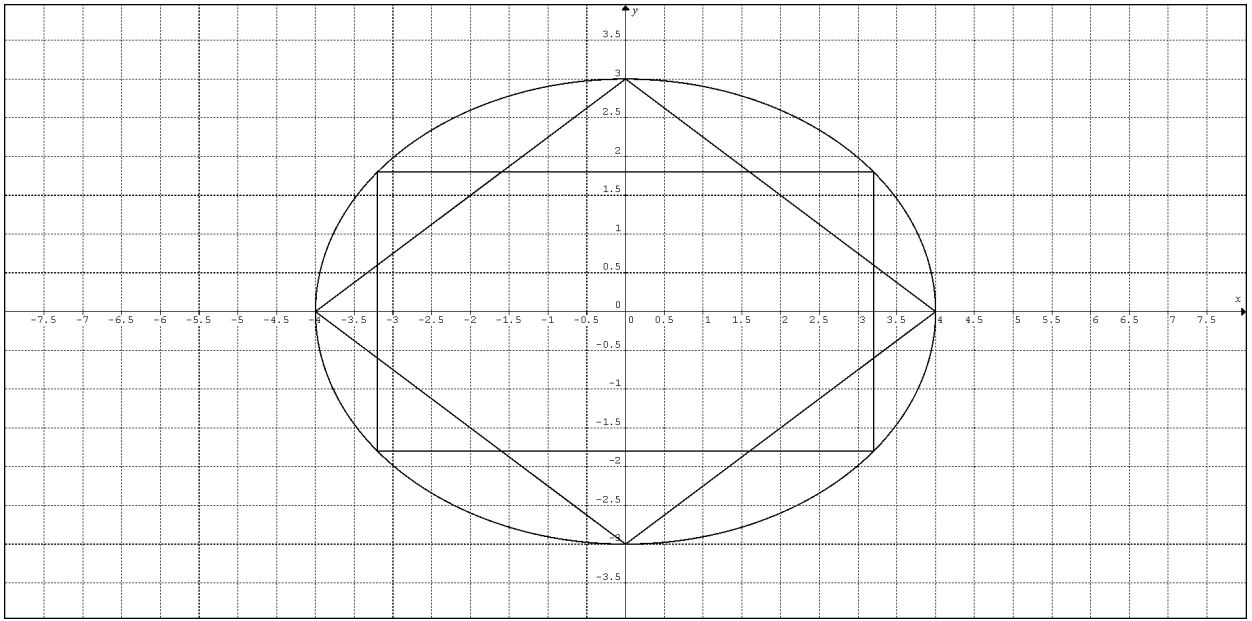

**Figure 1**

**Problem 1.** Is it true that the perimeter of any quadrilateral inscribed in the ellipse $\frac{x^2}{a^2}+\frac{y^2}{b^2}=1$ does not exceed $4\sqrt{a^2+b^2}$ ? If it is true, are there any quadrilaterals with perimeter $4\sqrt{a^2+b^2}$, except the rectangle of the greatest perimeter and the rhombus?

(c) Find the dimensions of the rectangle inscribed in the ellipse $\frac{x^2}{a^2}+\frac{y^2}{b^2}=1$ with the greatest ratio $R=\frac{S}{P^2}$, where $S$ is the area of the rectangle and $P$ is its perimeter. We have to maximize the ratio $\frac{xy}{(x+y)^2}$ subject to $\frac{x^2}{a^2}+\frac{y^2}{b^2}=1, x>0, y>0$. Equivalently, we have to minimize $\frac{(x+y)^2}{xy}=\frac{x}{y}+\frac{y}{x}+2$.

Using Lagrange multipliers, we write

$$\frac{1}{y}-\frac{y}{x^2}=\frac{2\lambda x}{a^2}$$

$$\frac{1}{x}-\frac{x}{y^2}=\frac{2\lambda y}{b^2}$$

From these equations and the equation of the ellipse follows that $\lambda = 0$.

Therefore, $x = y = \dfrac{ab}{\sqrt{a^2+b^2}}$ and $R_{\max} = \dfrac{1}{16}$.

**Problem 2.** Is it true that for any quadrilateral inscribed in the ellipse $\dfrac{x^2}{a^2}+\dfrac{y^2}{b^2}=1$ the ratio $R = \dfrac{S}{P^2}$ does not exceed $\dfrac{1}{16}$? If it is true, are there any quadrilaterals, inscribed in the ellipse, except the square, for which this value of $R$ is attained?

**Part 1B. The curves** $\dfrac{|x|^\alpha}{A^\alpha}+\dfrac{|y|^\alpha}{B^\alpha}=1, \alpha > 0$ (1)

(a) The same reasoning as in Part 1A (a) shows that the dimensions of the rectangle of the greatest area inscribed in the curve $\dfrac{|x|^\alpha}{A^\alpha}+\dfrac{|y|^\alpha}{B^\alpha}=1$ are defined by the equations $\dfrac{x^\alpha}{A^\alpha}=\dfrac{y^\alpha}{B^\alpha}=\dfrac{1}{2}$. Hence, $x = 2^{-1/\alpha}A, y = 2^{-1/\alpha}B$, and the maximal area is $4xy = 2^{2-2/\alpha}AB$.

**Remark 2.** It follows from the definition of the Beta function (see e.g. [Ar]) that the area of the region bounded by the curve $\dfrac{|x|^\alpha}{A^\alpha}+\dfrac{|y|^\alpha}{B^\alpha}=1$ is $4AB\dfrac{\mathbf{B}\left(\dfrac{1}{\alpha},\dfrac{1}{\alpha}+1\right)}{\alpha}$.

Thus, the ratio of the greatest possible area of a rectangle inscribed into the curve (1) and the area of the region bounded by this curve is $\mathbf{R}_{\max}(\alpha) = \dfrac{\alpha}{2^{2/\alpha}\mathbf{B}\left(\dfrac{1}{\alpha},\dfrac{1}{\alpha}+1\right)}$.

We leave it to the reader to verify that $\mathbf{R}_{\max}(\alpha)$ is a strictly decreasing function of $\alpha$ and that $\lim\limits_{\alpha \to 0+}\mathbf{R}_{\max}(\alpha) = 0$ and $\lim\limits_{\alpha \to \infty}\mathbf{R}_{\max}(\alpha) = 1$.

(b) The problem of inscribing the rectangle of the greatest (or smallest) possible perimeter into the curve (1) can be stated as optimizing the sum $x+y$ subject to $\frac{x^\alpha}{A^\alpha} + \frac{y^\alpha}{B^\alpha} = 1, x > 0, y > 0$. Using the Lagrange multipliers, we obtain

$\frac{\lambda \alpha x^{\alpha-1}}{A^\alpha} = \frac{\lambda \alpha y^{\alpha-1}}{B^\alpha} = 1$. Combining these equations with (1) we obtain

$$x = \frac{A^{\alpha/(\alpha-1)}}{\left(A^{\alpha/(\alpha-1)} + B^{\alpha/(\alpha-1)}\right)^{1/\alpha}}, y = \frac{B^{\alpha/(\alpha-1)}}{\left(A^{\alpha/(\alpha-1)} + B^{\alpha/(\alpha-1)}\right)^{1/\alpha}},$$ (2)

$$P = 4(x+y) = 4\left(A^{\alpha/(\alpha-1)} + B^{\alpha/(\alpha-1)}\right)^{(\alpha-1)/\alpha}$$

To find out whether these values of $x$ and $y$ provide the maximum or the minimum of $x + y$ we look at the bordered Hessian. We refer the reader to [Ch, pages 383, 384] for the definition of bordered Hessian. In our case the bordered Hessian is

$$H = \begin{vmatrix} 0 & \frac{\alpha x^{\alpha-1}}{A^\alpha} & \frac{\alpha y^{\alpha-1}}{B^\alpha} \\ \frac{\alpha x^{\alpha-1}}{A^\alpha} & \frac{\alpha(1-\alpha)x^{\alpha-2}}{A^\alpha} & 0 \\ \frac{\alpha y^{\alpha-1}}{B^\alpha} & 0 & \frac{\alpha(1-\alpha)y^{\alpha-2}}{B^\alpha} \end{vmatrix} = \frac{\alpha^3(\alpha-1)(xy)^{\alpha-2}(x^\alpha B^\alpha + y^\alpha A^\alpha)}{A^{2\alpha} B^{2\alpha}}$$

Therefore, if $\alpha > 1$ then $H > 0$ and formulas (2) define the inscribed rectangle of the greatest perimeter, while if $\alpha < 1$ then $H < 0$ and these formulas define the inscribed rectangle of the smallest perimeter.

**Remark 3.** Instead of using the bordered Hessian we can notice that when $x = 0$ or $y = 0$ the perimeter of the (degenerated) rectangle is $4B$, respectively $4A$, and that

$$4\left(A^{\alpha/(\alpha-1)} + B^{\alpha/(\alpha-1)}\right)^{(\alpha-1)/\alpha} > 4\max(A,B) \text{ if } \alpha > 1 \text{ and}$$

$$4\left(A^{\alpha/(\alpha-1)} + B^{\alpha/(\alpha-1)}\right)^{(\alpha-1)/\alpha} < 4\min(A,B) \text{ if } \alpha < 1.$$

**Remark 4.** It is trivial to notice that if $\alpha = 1$ and $A \neq B$ then neither the inscribed rectangle of the greatest perimeter, nor the one of the smallest perimeter, exist. While if $A = B$, all the inscribed rectangles have the same perimeter.

Figure 2 below shows the rectangle of the greatest perimeter inscribed into the curve $|x|^3 + \frac{|y|^3}{8} = 1$. The coordinates of the vertex of the rectangle in the first quadrant are $x = \frac{1}{\sqrt[3]{1+2\sqrt{2}}} \approx 0.6392340079$, $y = \frac{2\sqrt{2}}{\sqrt[3]{1+2\sqrt{2}}} \approx 1.808026807$.

The maximum perimeter is $4(x+y) = 4(1+2\sqrt{2})^{2/3} \approx 9.789043256$

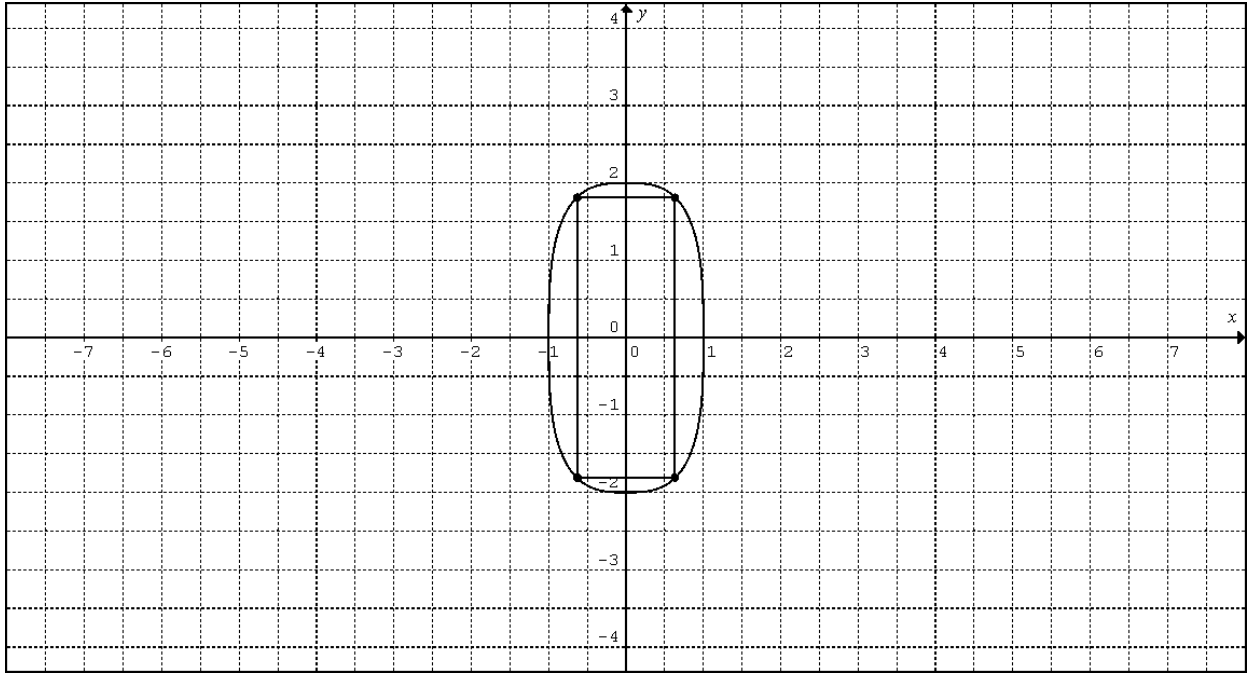

Figure 2

On Figure 3 we see the rectangle of the smallest perimeter inscribed in the curve $\sqrt[3]{x} + \frac{\sqrt[3]{y}}{2} = 1$. The coordinates of the vertex of the rectangle in the first quadrant are $x = \frac{1}{(1+1/2\sqrt{2})^3} \approx 0.4032494905$, $y = \frac{x}{2\sqrt{2}} \approx 0.1425702245$.

The perimeter of this rectangle is $P_{\min} = \frac{1}{(1+1/2\sqrt{2})^2} \approx 2.183278859$.

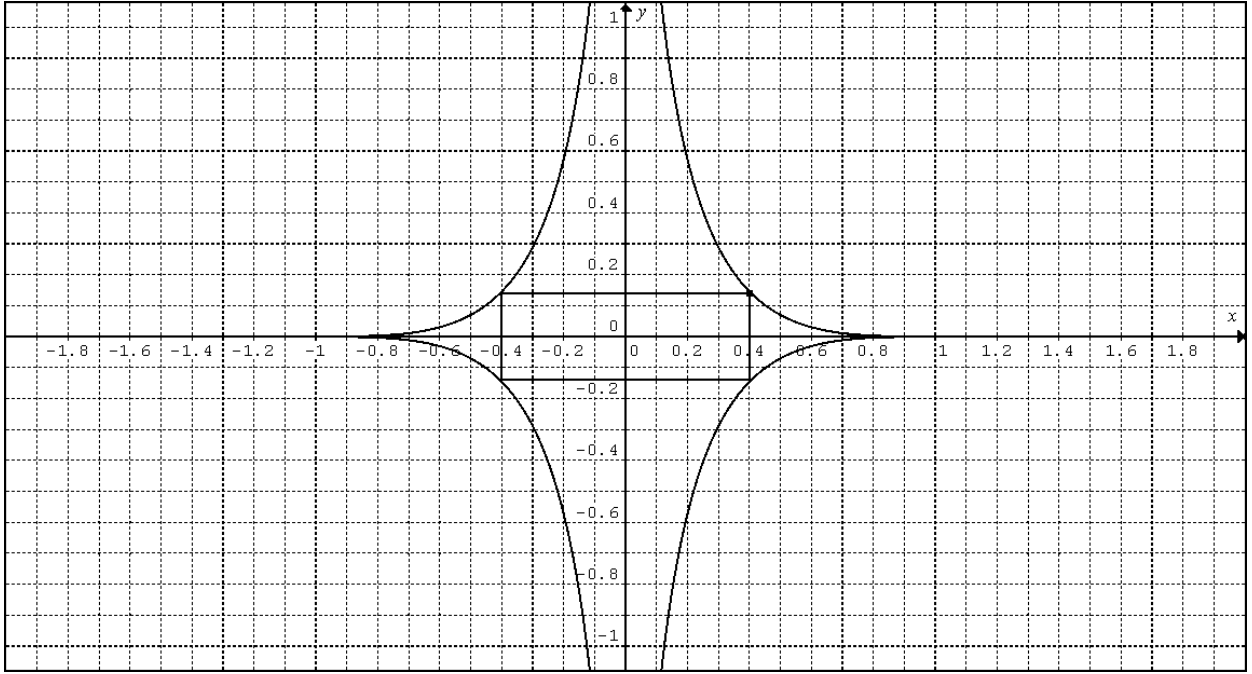

**Figure 3**

(c) We consider the problem of inscribing in the curve (1) the rectangle with the greatest ratio $R = \dfrac{S}{P^2} = \dfrac{xy}{4(x+y)^2}$. Equivalently we have to minimize

$\dfrac{(x+y)^2}{xy} = \dfrac{x}{y} + \dfrac{y}{x} + 2$ subject to $\dfrac{x^\alpha}{A^\alpha} + \dfrac{y^\alpha}{B^\alpha} = 1, x > 0, y > 0$. The same reasoning as in

**Part 1A** (c) shows that the maximum value of $R$, $R_{max} = \dfrac{1}{16}$ is attained when

$x = y$.

**Part 1C.** The curves $\dfrac{|x|^\alpha}{A^\alpha} + \dfrac{|y|^\beta}{B^\alpha} = 1, \alpha, \beta > 0$ (2).

(a) To find the dimensions of the rectangle of the greatest area inscribed in a curve (2) we use the method of Lagrange multipliers:

$$y = \dfrac{\alpha \lambda x^{\alpha-1}}{A^\alpha}, x = \dfrac{\beta \lambda y^{\beta-1}}{B^\beta} \quad (3)$$

From (3) follows $\dfrac{\alpha \lambda x^{\alpha}}{A^{\alpha}} = \dfrac{\beta \lambda y^{\beta}}{B^{\beta}} = xy$ and therefore,

$$\dfrac{y^{\beta}}{B^{\beta}} = \dfrac{x^{\alpha}}{A^{\alpha}} \cdot \dfrac{\alpha}{\beta} \quad (4)$$

Combining (2) and (4) we obtain

$$x = A\left[\dfrac{\beta}{(\alpha+\beta)}\right]^{1/\alpha}, \quad y = B\left[\dfrac{\alpha}{(\alpha+\beta)}\right]^{1/\beta}$$

$$S_{\max} = 4AB \dfrac{\alpha^{1/\beta} \beta^{1/\alpha}}{(\alpha+\beta)^{1/\alpha+1/\beta}}.$$

(b) Optimize the sum $x+y$ subject to $\dfrac{x^{\alpha}}{A^{\alpha}} + \dfrac{y^{\beta}}{B^{\beta}} = 1, x > 0, y > 0$.

**Proposition 1** (1) If $\alpha > 1$ and $\beta > 1$, then the rectangle of the greatest perimeter exists and is unique. The rectangle of the smallest perimeter does not exist.

(2) If $\alpha < 1$ and $\beta < 1$, then the rectangle of the smallest perimeter exists and is unique. The rectangle of the greatest perimeter does not exist.

(3) If $\alpha > 1$ and $\beta < 1$ (or $\alpha < 1$, and $\beta > 1$ ), then the rectangle of the greatest perimeter and the rectangle of the smallest perimeter both exist and are unique.

*Proof.* Let $x = Au, y = Bv$. Then $x + y = Au + Bv = Au + B(1-u^{\alpha})^{\delta}$, where $\delta = 1/\beta$. Let $F(u) = Au + B(1-u^{\alpha})^{\delta}$. Then

$$F'(u) = A - B\alpha\delta u^{\alpha-1}(1-u^{\alpha})^{\delta-1},$$
$$F''(u) = -B\alpha\delta u^{\alpha-2}(1-u^{\alpha})^{\delta-2}[(\alpha-1)(1-u^{\alpha}) + (1-\delta)\alpha u^{\alpha}]$$

We have to consider three cases.

1. Let $\alpha > 1, \beta > 1$. Because in this case $F'(0) = F'(1) = A > 0$ and $F'' < 0$ on $[0,1]$ the function $F$ has the unique maximum in $(0,1)$. Thus, the rectangle of the greatest perimeter exists and is unique. The rectangle of the smallest perimeter does not exist.

2. Let $\alpha < 1, \beta < 1$. In this case $\lim\limits_{u \to 0+} F'(u) = \lim\limits_{u \to 1-} F'(u) = -\infty$ and $F'' > 0$ on $[0,1]$. Therefore, the function $F$ has the unique minimum in (0,1). Thus, the rectangle of the smallest perimeter exists and is unique. The rectangle of the greatest perimeter does not exist.

3. Let $\alpha > 1, \beta < 1$ (or $\alpha < 1, \beta > 1$). In this case $F'(0) > 0$ and $\lim\limits_{u \to 1-} F'(u) = -\infty$. Moreover, in this case the equation $F''(u) = 0$ has exactly one solution on (0,1): $u = \left(\dfrac{\alpha - 1}{\alpha \delta - 1}\right)^{1/\alpha}$. Therefore $F$ has exactly one maximum and exactly one minimum on (0,1), and both the rectangle of the greatest perimeter and of the smallest perimeter exist and are unique. □

**Example 1.** We consider the curve $\dfrac{|x|^3}{27} + \dfrac{\sqrt{|y|}}{2} = 1$. We have to optimize $x + y$ subject to $\dfrac{x^3}{27} + \dfrac{\sqrt{y}}{2} = 1, x > 0, y > 0$. The method of Lagrange multipliers provides $1 = \lambda \dfrac{x^2}{9} = \lambda \dfrac{1}{4\sqrt{y}}$. Hence, $\sqrt{y} = \dfrac{9}{x^2}$, and from the equation of the curve we get

$$\dfrac{8}{27}x^5 - 8x^2 + 9 = 0. \quad (5)$$

Equation (5) has two positive solutions: $x_1 \approx 1.086810566$, $x_2 \approx 2.855105222$. The corresponding values of $y$ are $y_1 \approx 3.628687998$, $y_2 \approx 0.07618624391$. Figure 4 below shows the curve $\dfrac{|x|^3}{27} + \dfrac{\sqrt{|y|}}{2} = 1$ and the rectangles of the greatest and the smallest perimeter inscribed into it.

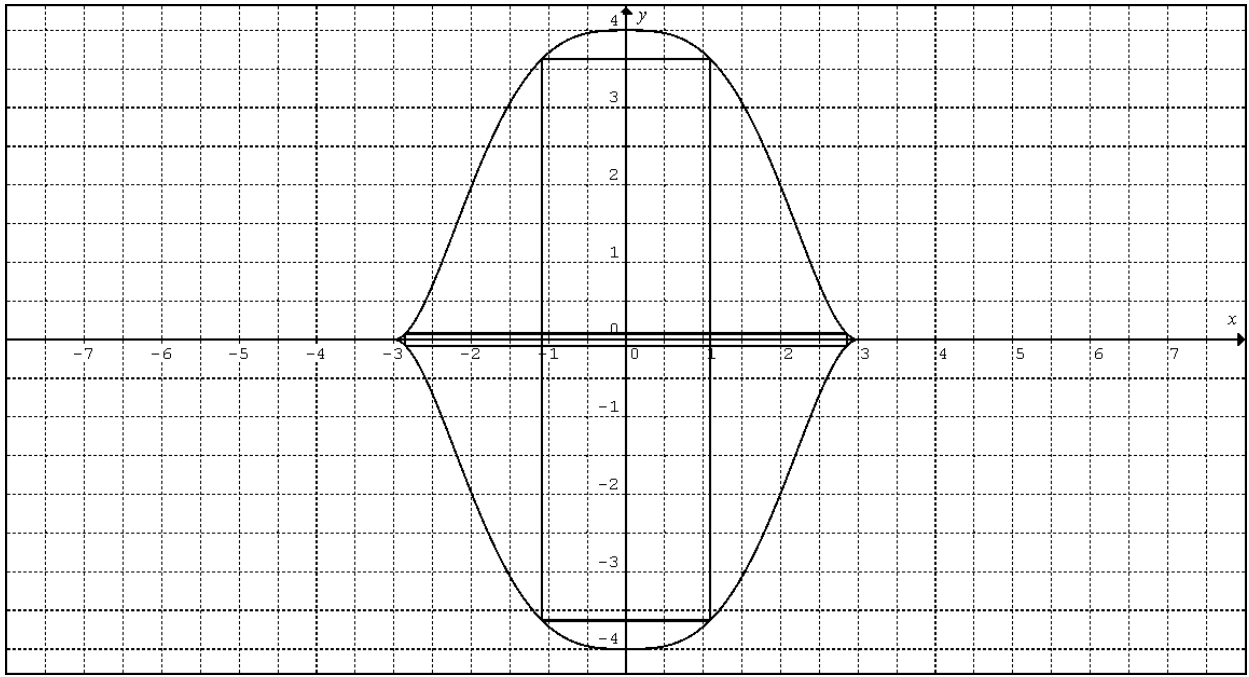

Figure 4

(c) Minimize $\dfrac{xy}{(x+y)^2}$ subject to $\dfrac{x^\alpha}{A^\alpha}+\dfrac{y^\beta}{B^\beta}=1$, $x>0, y>0$. The method of Lagrange multipliers provides the equations:

$$\frac{1}{y}-\frac{y}{x^2}=\frac{\lambda\alpha x^{\alpha-1}}{A^\alpha}$$

$$\frac{1}{x}-\frac{x}{y^2}=\frac{\lambda\beta y^{\beta-1}}{B^\beta}$$

From these equations follows that $\dfrac{\alpha\lambda x^\alpha}{A^\alpha}+\dfrac{\beta\lambda y^\beta}{B^\beta}=0$. Therefore, $\lambda=0$ and $x=y$.

**Part 1D. Regions bounded by closed curves symmetric about the $x$ and $y$-axes.**

In this part we assume that $A$ and $B$ are positive real numbers and that $f$ is a function continuous on $[0, A]$, strictly decreasing on this interval, and such that $f(0) = B, f(A) = 0$. We consider the curve defined by the equation:

$$|y| = f(|x|) \quad (6)$$

(a) The existence of a rectangle of the greatest area inscribed in a curve (6) follows from the standard argument: the function $xf(x)$ takes its greatest value in the interval $(0, A)$. Only the continuity of the decreasing function $f$ is required. More interesting is the question under what $f$ conditions the rectangle of the greatest area is unique. The following proposition is a trivial consequence of Rolle's theorem.

**Proposition 2.** Assume that the function $f$ is twice differentiable on $(0, A)$ and that at any point $t \in (0, A)$ we have $2f'(t) + tf''(t) \neq 0$. Then the rectangle of the greatest area inscribed in the curve (6) is unique.

**Corollary 1.** If the function $f$ is convex down, then the rectangle of the greatest area inscribed in the curve (6) is unique.

**Remark 5.** The sufficient condition in Proposition 2 is not necessary. Indeed, as we have seen in Part 1C (a), the rectangle of the greatest area inscribed in curve (2) is unique, but as the proof of Proposition 1 shows, in the case when $\min(\alpha, \beta) < 1$ the function $y = y(x)$ either is convex up or has exactly one inflection point in the interval $(0, A)$.

Nevertheless, as the next example shows, we cannot completely dispense with the condition in Proposition 1.

**Example 2.** There is a strictly decreasing on $[0, A]$, continuously differentiable on $(0, A)$ and strictly convex up function $f$ such that there are exactly two rectangles of the greatest area inscribed in the curve $|y| = f(|x|)$. We define the function $f$ on the interval $[0,4]$ as follows:

$$f(x) = \begin{cases} 4(1-\sqrt{x})^2, & 0 \le x \le \dfrac{9}{25} \\ P(x) = -\dfrac{15625}{1176}x^3 + \dfrac{14125}{588}x^2 - \dfrac{17401}{1176}x + \dfrac{851}{245}, & \dfrac{9}{25} < x \le \dfrac{16}{25} \\ \dfrac{1}{4}(2-\sqrt{x})^2, & \dfrac{16}{25} < x \le 4 \end{cases}$$

Figure 5 below shows the curve $|y| = f(|x|)$ and two inscribed in this curve rectangles of the greatest area. The vertices of these rectangles in the first quadrant are $\left(\dfrac{1}{4}, 1\right)$ and $\left(1, \dfrac{1}{4}\right)$, respectively.

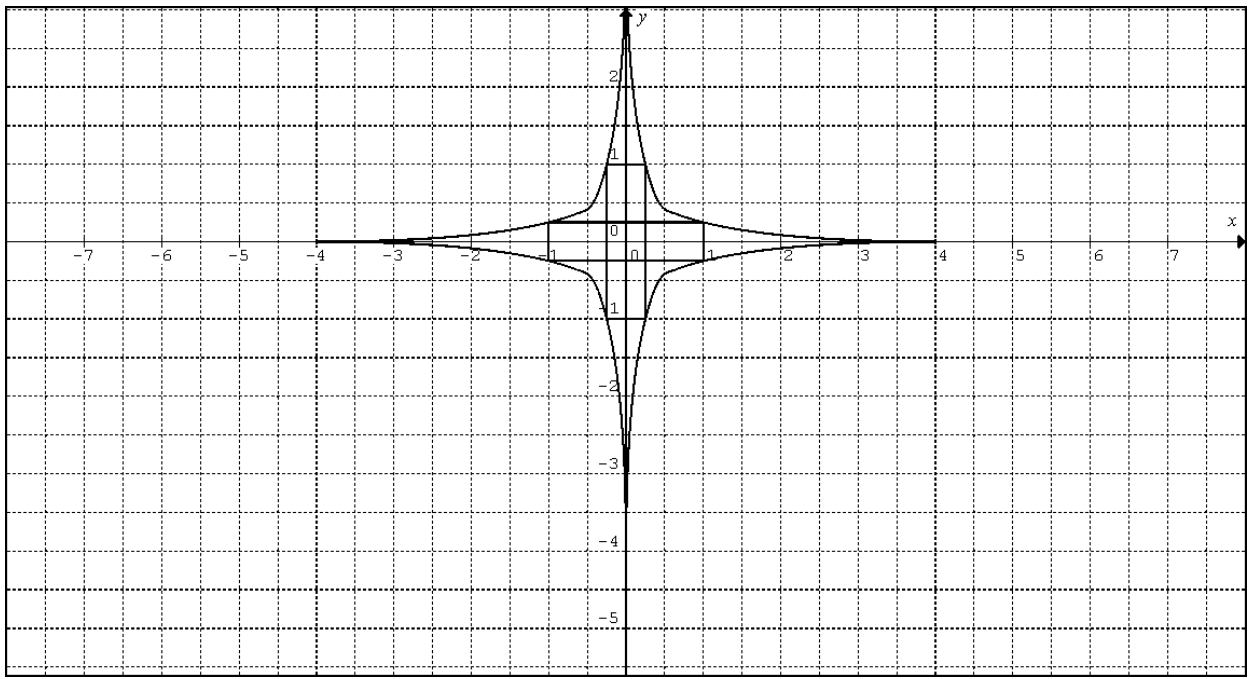

Figure 5

**Remark 6.** At the price of increasing the degree of polynomial $P$ we can for any positive integer $n$ find a function $f$ with the properties described in Example 2 and continuously differentiable $n$ times on $(0,4)$.

Clearly, if the decreasing function $f$ is real-analytic on $(0,A)$, then the problem of finding a rectangle of the greatest volume inscribed in the curve $|y| = f(|x|)$ has a finite number of solutions. In connection with it we pose the following problem.

**Problem 3.** Let $P$ be a polynomial function of degree $n$ decreasing on an interval $[0,A]$ and such that $P(A) = 0$. What is the maximum number of rectangles of the greatest area inscribed in the curve $|y| = P(|x|)$.

On the other hand, the situation is completely different if instead of real-analytic functions we consider functions continuously differentiable on $(0,A)$. We assume that the reader is familiar with some standard facts concerning the topology of the real line.

**Proposition 3.** Let $E$ be a nonempty closed bounded subset of $(0,\infty)$. There are a positive real number $A$ and a function $f$, such that:

(*) $E \subset (0, A)$,

(**) The function $f$ is continuous and strictly decreasing on $[0, A]$, moreover $f(A) = 0$

(***) The function $f$ is continuously differentiable and strictly convex up on $(0, A)$,

(****) The rectangle with the vertex $(x, f(x))$ inscribed in the curve $|y| = f(|x|)$ has the greatest possible area if and only if $x \in E$.

*Proof.* Let $m = \min\{x : x \in E\}$, and $M = \max\{x : x \in E\}$. The open set $[m, M] \setminus E$ is the union of at most countable family of pairwise disjoint open intervals, $E = \bigcup_{n=1}^{\infty} (a_n, b_n)$. We order these intervals in such a way that $b_{n+1} - a_{n+1} \leq b_n - a_n$. We include the case when the family of intervals $(a_n, b_n)$ is finite, or $E = [m, M]$ by assuming that for some $p \geq 1$ the interval $(a_p, b_p)$ is empty. We define the function $f$ as follows.

$$f(x) = \begin{cases} \dfrac{1}{x}, & \text{if } x \in E \\[4pt] \dfrac{1}{x} - \varepsilon_n (x - a_n)^2 (x - b_n)^2, & \text{if } x \in (a_n, b_n) \\[4pt] \dfrac{1}{m} - \dfrac{1}{m^2}(x - m) + \dfrac{1}{m^3}(x - m)^2, & \text{if } x \in [0, m) \\[4pt] \dfrac{1}{M} - \dfrac{1}{M^2}(x - M) + \varepsilon (x - M)^2, & \text{if } x \in (M, A], \text{ where } \varepsilon \leq \dfrac{1}{4M^3} \\[4pt] \text{and } A = \dfrac{2\varepsilon M^3 + 1 + \sqrt{1 - 4\varepsilon M^3}}{2\varepsilon M^2} \end{cases} \qquad (7)$$

The positive numbers $\varepsilon_n$ in (7) clearly can be chosen in such a way that the function $f$ is decreasing and strictly convex up. Moreover, if the family of intervals $(a_n, b_n)$ is infinite, we chose $\varepsilon_n$ in such a way that $\varepsilon_n \downarrow 0$, ensuring that $f$ is continuously differentiable on $(0, 2M)$.

**Remark 7.** It is not difficult to modify the definition of the function $f$ in such a way that $f$ becomes continuously differentiable $n$ times where $n$ is an arbitrary positive integer. Moreover, if the family of intervals $(a_n, b_n)$ is finite we can modify the definition of $f$ in such a way that $f$ becomes differentiable infinitely many times on $(0, A)$. Indeed, we put,

$$f(x) = \frac{1}{x} + \varepsilon_n \exp\left(-\frac{1}{(x - a_n)^2}\right) \exp\left(-\frac{1}{(x - b_n)^2}\right) \quad x \in (a_n, b_n).$$

The definition of $f$ on $[0, m]$ and on $[M, A]$ should be modified as well.

Remark 7 gives rise to our last problem in this subsection.

**Problem 4.** Let $E$ be a closed bounded subset of $(0, \infty)$. What conditions on $E$ are necessary and/or sufficient for existence of a function $f$ with properties listed in Proposition 3 and differentiable infinitely many times on $(0, A)$?

(b) In regard to inscribing a rectangle of the largest or the smallest perimeter in a curve (6) we state the following trivial proposition.

**Proposition 4.** Let $A > 0$ and $f$ be a strictly decreasing continuous function on $[0, A]$ such that $f(A) = 0$. Let $C$ be the closed curve defined by the equation $|y| = f(|x|)$.

(I) Assume that $f$ is twice differentiable and convex down on $(0, A)$. Also assume that $\lim_{x \to 0+} f'(x) > -1$ and $\lim_{x \to A-} f'(x) < -1$, then the rectangle of the greatest perimeter inscribed in $C$ exists and is unique.

(II) Similarly, if $f$ is twice differentiable and convex up on $(0, A)$ and $\lim_{x \to 0+} f'(x) < -1$ and $\lim_{x \to A-} f'(x) > -1$, then the rectangle of the smallest perimeter inscribed in $C$ exists and is unique.

(III) Let $E$ be a closed subset of $[0, A]$ and $n$ is a positive integer. There is a $n$ times continuously differentiable on $(0, A)$, strictly decreasing function $f$ such that $f(0) = A, f(A) = 0$, and the rectangle inscribed in the curve $|y| = f(|x|)$ with the vertex $(x, f(x))$ has the greatest possible perimeter if and only if $x \in E \cap (0, A)$.

**Problem 5.** Let $A$ and $B$ be positive real numbers. Describe the classes **A** and **B** of closed subsets of $[0, A]$ defined as follows.

$E \in \mathbf{A}$ if and only if there is a strictly decreasing function $f$ such that $f(0) = B, f(A) = 0$, and the rectangle inscribed in the curve $|y| = f(|x|)$ with the vertex $(x, f(x))$ has the **greatest** possible perimeter if and only if $x \in E \cap (0, A)$.

$E \in \mathbf{B}$ if and only if there is a strictly decreasing function $f$ such that $f(0) = B, f(A) = 0$, and the rectangle inscribed in the curve $|y| = f(|x|)$ with the vertex $(x, f(x))$ has the **smallest** possible perimeter if and only if $x \in E \cap (0, A)$.

(c) In this subsection we consider the problem of inscribing in curve (6) a rectangle with the greatest ratio $\mathbf{R} = \dfrac{S}{P^2}$.

**Proposition 5.** Let $f$ be a continuous strictly decreasing function on $[0, A]$ such that $f(A) = 0$. Among the rectangles inscribed in the curve $|y| = f(|x|)$ the square has the greatest ratio $\mathbf{R}$.

*Proof.* Assume first that the function $f$ is continuously differentiable on $[0, A]$. We have to minimize $\dfrac{x}{y} + \dfrac{y}{x}$ subject to $y - f(x) = 0$. The method of Lagrange multipliers provides the following equations:

$$\frac{1}{y} - \frac{y}{x^2} = \lambda f'(x)$$

$$\frac{1}{x} - \frac{x}{y^2} = -\lambda$$

Therefore, $x^2 - y^2 = \lambda x^2 y f'(x) = \lambda x y^2$. Hence, $x = y$.

The general case follows from what we have proved because, if $f$ is an arbitrary strictly decreasing continuous function on $[0, A]$ such that $f(A) = 0$, then for any integer $n$ we can find a continuously differentiable and strictly decreasing on $[0, A]$ function $f_n$ such that $f_n(A) = 0$ and $\max_{x \in [0, A]} |f(x) - f_n(x)| \leq \dfrac{1}{n}$ □

**Remark 8.** The reader not familiar with the Weierstrass approximation theorem can use the following simple argument in the proof of Proposition 5. Divide $[0, A]$ into small intervals, on each of these intervals change $f$ to the linear function with the same values at the ends of the interval and smooth the curve at the angles by using appropriate cubic polynomials.

# Part 2. Optimization problems in three dimensions.

## Part 2A. Inscribing a rectangular parallelepiped with the greatest volume in a surface symmetric about coordinate planes.

Let $\alpha > 0$. The problem of inscribing the rectangular parallelepiped of the greatest volume in the surface $\dfrac{|x|^\alpha}{A^\alpha} + \dfrac{|y|^\alpha}{B^\alpha} + \dfrac{|z|^\alpha}{C^\alpha} = 1$ is trivial. Because the product of three positive numbers with the fixed sum takes the greatest value if and only if the numbers are equal, we have $\dfrac{|x|^\alpha}{A^\alpha} = \dfrac{|y|^\alpha}{B^\alpha} = \dfrac{|z|^\alpha}{C^\alpha} = \dfrac{1}{3}$. Therefore,

$x = \dfrac{A}{3^{1/\alpha}}, y = \dfrac{B}{3^{1/\alpha}}, z = \dfrac{C}{3^{1/\alpha}}$, and $V_{max} = \dfrac{8ABC}{3^{3/\alpha}}$.

**Example 3.** Let $\alpha = 3$. Then $V_{max} = \dfrac{8ABC}{3}$. The volume of the solid bounded by the surface $\dfrac{|x|^3}{A^3} + \dfrac{|y|^3}{B^3} + \dfrac{|z|^3}{C^3} = 1$ is $V_{solid} = 8ABC \dfrac{8\pi^3 \sqrt{3}}{243\Gamma\left(\dfrac{2}{3}\right)^3}$. The ratio $\dfrac{V_{max}}{V_{solid}}$ is

$\dfrac{81\Gamma\left(\dfrac{2}{3}\right)^3}{8\pi^3 \sqrt{3}} \approx 0.4681168362$.

Next, we consider the region in $\mathbf{R}^3$ bounded by the surface

$$\dfrac{|x|^\alpha}{A^\alpha} + \dfrac{|y|^\beta}{B^\beta} + \dfrac{|z|^\gamma}{C^\gamma} = 1, \quad A, B, C, \alpha, \beta, \gamma > 0. \tag{8}$$

We have to maximize the product $xyz$ subject to $\dfrac{x^\alpha}{A^\alpha} + \dfrac{y^\beta}{B^\beta} + \dfrac{z^\gamma}{C^\gamma} = 1, \quad x, y, z \geq 0$.

The method of Lagrange multipliers provides the following equations:

$$yz = \dfrac{\alpha\lambda x^{\alpha-1}}{A^\alpha}, xz = \dfrac{\beta\lambda y^{\beta-1}}{B^\beta}, xy = \dfrac{\gamma\lambda z^{\gamma-1}}{C^\gamma}$$

From these equations we obtain that $\dfrac{\alpha x^\alpha}{A^\alpha} = \dfrac{\beta y^\beta}{B^\beta} = \dfrac{\gamma z^\gamma}{C^\gamma}$. Combining these equations with the equation (8) of the surface we get

$$x = A\left(\frac{\beta\gamma}{\alpha\beta + \alpha\gamma + \beta\gamma}\right)^{1/\alpha}, y = B\left(\frac{\alpha\gamma}{\alpha\beta + \alpha\gamma + \beta\gamma}\right)^{1/\beta}, z = C\left(\frac{\alpha\beta}{\alpha\beta + \alpha\gamma + \beta\gamma}\right)^{1/\gamma}.$$

The inscribed rectangular parallelepiped of maximum volume is unique, and its volume is
$$V_{max} = 8xyz = 8ABC \frac{\alpha^{(1\beta+1/\gamma)}\beta^{(1\alpha+1/\gamma)}\gamma^{(1\alpha+1/\beta)}}{(\alpha\beta + \alpha\gamma + \beta\gamma)^{(1/\alpha+1/\beta+1/\gamma)}}.$$

Assume that **S** is a compact connected surface in $\mathbf{R}^3$ symmetric about coordinate planes. Assume also that every ray starting at the origin intersects **S** exactly one time. The existence of a rectangular parallelepiped of the greatest volume inscribed in **S** with faces parallel to coordinate planes follows from a standard compactness argument. The following simple proposition singles out a special case when the rectangular parallelepiped of the greatest volume is unique.

**Proposition 6.** Let $f, g, h$ be functions continuous on $[0, A]$, twice differentiable on $(0, A)$ and convex up. Also assume that $f(0) = g(0) = h(0) = 0$. Let **S** be the surface in $\mathbf{R}^3$ defined by the equation,

$$f(|x|) + g(|y|) + h(|z|) = 1 \qquad (9)$$

Then the rectangular parallelepiped of the greatest volume inscribed in **S** (with faces parallel to coordinate planes) is unique.

*Proof.* We have to maximize the product $xyz$ subject to the constraints $f(x) + g(y) + h(z) = 1$, $x > 0, y > 0, z > 0$. The method of Lagrange multipliers provides the following equations:

$$\begin{aligned} yz &= \lambda f'(x) \\ xz &= \lambda g'(y) \\ xy &= \lambda h'(z) \end{aligned} \qquad (10)$$

Assume, contrary to our claim, that equations (10) are satisfied at two distinct points in the first octant: $(x_0, y_0, z_0)$ and $(x_1, y_1, z_1)$. We can assume without loss of generality that $z_1 > z_0$. From (10) we obtain:

$$\begin{aligned} y_0 g'(y_0) &= x_0 f'(x_0) & y_1 g'(y_1) &= x_1 f'(x_1) & (a) \\ z_0 h'(z_0) &= y_0 g'(y_0) & z_1 h'(z_1) &= y_1 g'(y_1) & (b) \quad (11) \\ z_0 h'(z_0) &= x_0 f'(x_0) & z_1 h'(z_1) &= x_1 f'(x_1) & (c) \end{aligned}$$

In virtue of our assumptions the functions $xf'(x), yg'(y),$ and $zh'(z)$ are strictly increasing on $[0, A]$. Indeed, $(xf'(x))' = f'(x) + xf''(x) > 0$. It follows from $z_1 > z_0$ and (11c) that $x_1 > x_0$. Therefore, from (11a) we have $y_1 > y_0$, But these inequalities and the fact that the functions $f, g, h$ are strictly increasing clearly contradict the equation (9) □

If we assume the conditions of Proposition 6, and assume that the second derivatives $f'', g'', h''$ are strictly positive on $(0, A)$ then the Gaussian curvature $K$ of the surface $S$ is positive. Indeed, (see [Sp])

$$K(x, y, z) = \frac{f''(x)g''(y)(h'(z))^2 + f''(x)h''(z)(g'(y))^2 + g''(y)h''(z)(f'(x))^2}{(f'(x)^2 + g'(y)^2 + h'(z)^2)^2}.$$

Therefore, the solid bounded by $S$ is convex and we pose the following problem.

**Problem 6.** Let $S$ be a compact closed smooth surface in $\mathbf{R}^3$ symmetric about the coordinate planes. Assume that the solid bounded by $S$ is strictly convex (the intersection of every plane tangent to $S$ with $S$ is a singleton). Is it true that the rectangular parallelepiped (with faces parallel to coordinate planes) of the greatest volume inscribed in $S$ is unique?

**Part 2B. Inscribing a rectangular parallelepiped of the greatest surface area in a surface symmetric about coordinate planes.**

This problem is considerably more complicated than the problem described in Part 2A. We will start by describing its complete solution in the case of an ellipsoid.

We look at the following optimization problem:

Maximize $xy + xz + yz$ subject to $\dfrac{x^2}{a^2} + \dfrac{y^2}{b^2} + \dfrac{z^2}{c^2} = 1, x > 0, y > 0, z > 0.$

The Lagrange equations are:

$$\begin{bmatrix} y+z = \dfrac{2\lambda x}{a^2} \\ x+z = \dfrac{2\lambda y}{b^2} \\ x+y = \dfrac{2\lambda z}{c^2} \end{bmatrix} \quad (12)$$

The system (12) of linear equations has a nontrivial solution if and only if the determinant:

$$\begin{vmatrix} \dfrac{2\lambda}{a^2} & -1 & -1 \\ -1 & \dfrac{2\lambda}{b^2} & -1 \\ -1 & -1 & \dfrac{2\lambda}{c^2} \end{vmatrix} = \dfrac{2(a^2b^2c^2 + \lambda(a^2b^2 + a^2c^2 + b^2c^2) - 4\lambda^3)}{a^2b^2c^2} = 0. \quad (13)$$

The equation $P(\lambda) = 4\lambda^3 - \lambda(a^2b^2 + a^2c^2 + b^2c^2) - a^2b^2c^2 = 0$ has exactly one positive solution $\Lambda$. We will prove that maximum surface area of the inscribed parallelepiped is $4\Lambda$. First notice that the system

$$\begin{bmatrix} y+z = \dfrac{2\Lambda x}{a^2} \\ x+z = \dfrac{2\Lambda y}{b^2} \\ x+y = \dfrac{2\Lambda z}{c^2} \end{bmatrix} \quad (14)$$

has the only one solution $(x_0, y_0, z_0)$ that satisfies the equation $\dfrac{x_0^2}{a^2} + \dfrac{y_0^2}{b^2} + \dfrac{z_0^2}{c^2} = 1$. Next notice that multiplying both parts of equations (14) by $x_0, y_0,$ and $z_0$, respectively, and adding the obtained equations we get $x_0 y_0 + x_0 z_0 + y_0 z_0 = \Lambda$. It remains to prove that the maximum of $xy + xz + yz$ cannot be attained on the boundary of the surface $\dfrac{x^2}{a^2} + \dfrac{y^2}{b^2} + \dfrac{z^2}{c^2} = 1$, $x, y, z \geq 0$. Let us assume that $z = 0$. The maximum value of $xy$, subject to the constraint $\dfrac{x^2}{a^2} + \dfrac{y^2}{b^2} = 1$, is $\dfrac{ab}{2}$. But $P\left(\dfrac{ab}{2}\right) = -\dfrac{acb^2(a+c)^2}{2} < 0$, and therefore $\dfrac{ab}{2} < \Lambda$. The cases $y = 0$ or $x = 0$ can be considered similarly.

In the case of an ellipsoid of revolution, assuming e.g. that $b=c$, we can solve the cubic equation (13) explicitly and obtain the following formula,

$$S_{max} = 4\Lambda = b(b+\sqrt{8a^2+b^2}).$$

We will now state and prove a more general result.

**Proposition 7.** The rectangular parallelepiped of the greatest surface area with faces parallel to coordinate planes inscribed in the surface:

$$\frac{|x|^\alpha}{A^\alpha} + \frac{|y|^\alpha}{B^\alpha} + \frac{|z|^\alpha}{C^\alpha} = 1,$$

where $\alpha > 1$, exists and is unique.

*Proof.* We have to maximize $xy + xz + yz$ subject to the constraint:

$$\frac{x^\alpha}{A^\alpha} + \frac{y^\alpha}{B^\alpha} + \frac{z^\alpha}{C^\alpha} = 1, x > 0, y > o, z > 0. \qquad (15)$$

We will divide the proof of Proposition 7 into two steps.

Step 1. The function $xy + xz + yz$ attains its global maximum at an inner point of the surface $\frac{x^\alpha}{A^\alpha} + \frac{y^\alpha}{B^\alpha} + \frac{z^\alpha}{C^\alpha} = 1, x > 0, y > o, z > 0$. Assume to the contrary that the global maximum is attained at some point $(x_0, y_0, z_0)$ on the boundary of surface (15). Without loss of generality, we can assume that $z_0 = 0$. Then (see Part1B (a)) we can assume that $x_0 = 2^{-1/\alpha} A, y_0 = 2^{-1/\alpha} B$. We fix an $\varepsilon > 0$ and consider the point $(x_1, y_1, z_1)$, where $x_1 = \frac{A}{2^{1/\alpha}} \left(1 - \frac{\varepsilon^\alpha}{C^\alpha}\right)^{1/\alpha}, y_1 = \frac{B}{2^{1/\alpha}} \left(1 - \frac{\varepsilon^\alpha}{C^\alpha}\right)^{1/\alpha}, z_1 = \varepsilon$. Then,

$$x_1 y_1 + x_1 z_1 + y_1 z_1 - x_0 y_0 = \frac{AB}{2^{2/\alpha}} \left[\left(1 - \frac{\varepsilon^\alpha}{C^\alpha}\right)^{2/\alpha} - 1\right] + \varepsilon \frac{A+B}{2^{1/\alpha}} \left(1 - \frac{\varepsilon^\alpha}{C^\alpha}\right)^{1/\alpha}.$$

The last expression is positive if $\varepsilon$ is a very small positive number. Indeed, by L'Hôpital's rule,

$$\lim_{\varepsilon \to 0+} \frac{\left[\left(1-\frac{\varepsilon^\alpha}{C^\alpha}\right)^{2/\alpha}-1\right]}{\varepsilon\left(1-\frac{\varepsilon^\alpha}{C^\alpha}\right)^{1/\alpha}} = \lim_{\varepsilon \to 0+} \frac{\left[\left(1-\frac{\varepsilon^\alpha}{C^\alpha}\right)^{2/\alpha}-1\right]}{\varepsilon} = \lim_{\varepsilon \to 0+} \left(1-\frac{\varepsilon^\alpha}{C^\alpha}\right)^{2/\alpha-1} \frac{-\alpha\varepsilon^{\alpha-1}}{C^\alpha} = 0.$$

Step 2. We prove that there is a unique point on the surface (15), $(X, Y, Z), X > 0, Y > 0, Z > 0$, at which is attained the global maximum of the function $xy + xz + yz$. Assume to the contrary that the global maximum is attained at two distinct points $(X_0, Y_0, Z_0)$ and $(X_1, Y_1, Z_1)$ on the surface (15). Then

$$X_0 Y_0 + X_0 Z_0 + Y_0 Z_0 = X_1 Y_1 + X_1 Z_1 + Y_1 Z_1 \qquad (16)$$

The method of Lagrange multipliers provides the following equations:

$$Y_0 + Z_0 = \frac{\alpha \lambda_0 X_0^{\alpha-1}}{A^\alpha}, \qquad Y_1 + Z_1 = \frac{\alpha \lambda_1 X_1^{\alpha-1}}{A^\alpha},$$

$$X_0 + Z_0 = \frac{\alpha \lambda_0 Y_0^{\alpha-1}}{B^\alpha}, \qquad X_1 + Z_1 = \frac{\alpha \lambda_1 Y_1^{\alpha-1}}{B^\alpha}, \qquad (17)$$

$$X_0 + Y_0 = \frac{\alpha \lambda_0 Z_0^{\alpha-1}}{C^\alpha}, \qquad X_1 + Y_1 = \frac{\alpha \lambda_1 Z_1^{\alpha-1}}{C^\alpha}.$$

Multiplying the equations in the left column by $X_0, Y_0$ and $Z_0$, respectively, adding them, and applying (15) we get $\alpha \lambda_0 = 2(X_0 Y_0 + X_0 Z_0 + Y_0 Z_0)$. Similarly, $\alpha \lambda_1 = 2(X_1 Y_1 + X_1 Z_1 + Y_1 Z_1)$. From (16) follows that $\lambda_1 = \lambda_0 = \lambda$ and equations (17) can be written as

$$Y_0 + Z_0 = \frac{\alpha \lambda X_0^{\alpha-1}}{A^\alpha}, \qquad Y_1 + Z_1 = \frac{\alpha \lambda X_1^{\alpha-1}}{A^\alpha}, \quad (a)$$

$$X_0 + Z_0 = \frac{\alpha \lambda Y_0^{\alpha-1}}{B^\alpha}, \qquad X_1 + Z_1 = \frac{\alpha \lambda Y_1^{\alpha-1}}{B^\alpha}, \quad (b) \qquad (18)$$

$$X_0 + Y_0 = \frac{\alpha \lambda Z_0^{\alpha-1}}{C^\alpha}, \qquad X_1 + Y_1 = \frac{\alpha \lambda Z_1^{\alpha-1}}{C^\alpha}. \quad (c)$$

We can assume without loss of generality that $Z_1 > Z_0$. Then it follows from (18c) that either $X_1 > X_0$ or $Y_1 > Y_0$. If e.g., $Y_1 > Y_0$ then it follows from (18a) that $X_1 > X_0$, a contradiction □

The condition $\alpha > 1$ is sufficient but in general not necessary for existence and/or uniqueness of the rectangular parallelepiped of the greatest surface area inscribed in the surface $\frac{|x|^{\alpha}}{A^{\alpha}}+\frac{|y|^{\alpha}}{B^{\alpha}}+\frac{|z|^{\alpha}}{C^{\alpha}}=1$. We were not able to find necessary and sufficient conditions on the parameter $\alpha$ (see Problem 7 below) but in our next proposition we provide such conditions in a special case when the surface is defined by the equation:

$$|x|^{\alpha}+|y|^{\alpha}+|z|^{\alpha}=1, \quad \alpha>0 \qquad (19)$$

**Proposition 8.** If $\alpha \geq \frac{\ln 3 - \ln 2}{2\ln 3}$ then the rectangular parallelepiped of the greatest surface area inscribed in the surface (19) exists and is unique. Moreover, it is the cube with sides equal to $2\left(\frac{1}{3}\right)^{1/\alpha}$.

If $0 < \alpha < \frac{\ln 3 - \ln 2}{2\ln 3} \approx 0.7381404932$ then the rectangular parallelepiped of the greatest surface area inscribed in the surface (19) does not exist.

*Proof.* The method of Lagrange multipliers provides the following equations.

$$y+z = \alpha\lambda x^{\alpha-1}, x+z = \alpha\lambda y^{\alpha-1}, x+y = \alpha\lambda z^{\alpha-1} \qquad (20)$$

We will prove that the only solution of system (20) is $x = y = z = (1/3)^{1/\alpha}$. We divide the proof into three parts.

1. $\alpha \geq 1$. From (20) we have $y - x = \lambda\alpha(x^{\alpha-1} - y^{\alpha-1})$, hence $y = x$. Similarly, $z = x$.

2. $\alpha = \frac{m}{n}, m, n \in N, m < n$.. Let $u = x^{1/n}, v = y^{1/n}, w = z^{1/n}$. Our problem then becomes to maximize $u^n v^n + u^n w^n + v^n w^n$ subject to $u^m + v^m + w^m = 1$. From (20) we obtain

$$nu^{n-1}v^n + nu^{n-1}w^n = \lambda m u^{m-1}$$
$$nv^{n-1}u^n + nv^{n-1}w^n = \lambda m v^{m-1}$$
$$nw^{n-1}u^n + nw^{n-1}v^n = \lambda m w^{m-1}$$

It follows that: $u^{n-m}(v^n + w^n) = v^{n-m}(u^n + w^n) = w^{n-m}(u^n + v^n)$

Assume, to the contrary that there is a solution distinct from $x = y = z = (1/3)^{1/\alpha}$. Without loss of generality, we can assume that $u > v \geq w$. Assume first that

$m \geq n - m$. The first equation can be written as $u^{n-m}v^{n-m}(u^m - v^m) = w^n(u^{n-m} - v^{n-m})$.

Using the elementary inequality $\dfrac{u^m - v^m}{u^{n-m} - v^{n-m}} \geq u^{2m-n}$ we obtain $\dfrac{u^{n-m}v^{n-m}(u^m - v^m)}{u^{n-m} - u^{n-m}} > w^n$, a contradiction. Assume now that $n - m > m$. From the second equation $v^{n-m}w^{n-m}(v^m - w^m) = u^n(v^{n-m} - w^{n-m})$ and therefore

$$u^n \frac{v^{n-m} - w^{n-m}}{v^m - w^m} \geq u^n v^{n-2m} > v^{n-m}w^{n-m},$$

again, a contradiction.

3. $0 < \alpha < 1$. Let $\alpha_n \to \alpha$ where $\alpha_n$ are rational numbers. From the previous step follows that the Lagrange equations for the problem: maximize $(fg)^{1/\alpha} + (fh)^{1/\alpha} + (gh)^{1/\alpha}$ subject to $f + g + h = 1, f, g, h \geq 0$ have the unique solution $f = g = h = 1/3$. It remains to no notice that

$(fg)^{1/\alpha_n} + (fh)^{1/\alpha_n} + (gh)^{1/\alpha_n} \to (fg)^{1/\alpha} + (fh)^{1/\alpha} + (gh)^{1/\alpha}$ uniformly on the surface $f + g + h = 1, f, g, h \geq 0$.

To finish the proof of proposition it remains to notice that the value of $xy + xz + yz$ at the point $(3^{-\alpha}, 3^{-\alpha}, 3^{-\alpha})$ is $3^{1-2\alpha}$. On the other hand, the greatest value of this function on the boundary of the region $x^\alpha + y^\alpha + \jmath^\alpha = 1, x > 0, y > 0, z > 0$ is $2^{-2/\alpha}$. The two expressions are equal when $3^{1-2/\alpha} = 2^{-2/\alpha}$, solving this equation we obtain $\alpha = \dfrac{2(\ln 3 - \ln 2)}{\ln 3} \approx 0.7381404932$ □

**Problem 7.** (1) Describe all the triples $(\alpha, \beta, \gamma)$ such that there exists the rectangular parallelepiped of the greatest surface area inscribed in the surface (8). For which triples such a rectangular parallelepiped exists and is unique?

(2). Let S be a compact closed smooth surface in $\mathbf{R}^3$ symmetric about coordinate planes. Assume that the solid bounded by S is strictly convex (the intersection of every plane tangent to S with S is a singleton). Is it true that the rectangular parallelepiped (with faces parallel to coordinate planes) of the greatest surface area inscribed in S exists and is unique?

We finish this subsection with two examples. The corresponding calculations were performed with the help of Maple and Mathematica.

**Example 4.** The parallelepiped of the greatest surface area inscribed into the surface $|x|^3 + \frac{|y|^3}{8} + \frac{|z|^3}{27} = 1$ exists and is unique. The coordinates of its vertex in the first octant are $x \approx 0.5842341046, y \approx 1.446885928, z \approx 2.250143153$. The surface area is $S_{max} = 4(xy + xz + yz) \approx 21.66252376$.

**Example 5.** The parallelepiped of the greatest surface area inscribed into the surface $x^2 + y^2 + |z|^3 = 1$ exists and is unique. The coordinates of its vertex in the first octant are $x \approx 0.611045, y \approx 0.611045, z \approx 0.632676$. The surface area is $S_{max} = 4(xy + xz + yz) \approx 4.586252$.

**Part 2C. Inscribing a rectangular parallelepiped with the greatest sum of the edges in a surface symmetric about coordinate planes.**

**Proposition 9**. Assume that $\alpha, \beta, \gamma > 1$. If the rectangular parallelepiped inscribed in the surface $\frac{|x|^\alpha}{A^\alpha} + \frac{|y|^\beta}{B^\beta} + \frac{|z|^\gamma}{C^\gamma} = 1$ with the greatest sum of the edges exists, then it is unique.

Proof. The method of Lagrange multipliers provides the following equations.

$$\frac{\alpha \lambda x^{\alpha-1}}{A^\alpha} = \frac{\beta \lambda y^{\beta-1}}{B^\beta} = \frac{\gamma \lambda z^{\gamma-1}}{C^\gamma} = 1. \qquad (21)$$

From these equations and the equation of the surface we get

$$\frac{x^\alpha}{A^\alpha} + \left(\frac{\alpha B^\beta}{\beta A^\alpha} x^{\alpha-1}\right)^{\beta/(\beta-1)} + \left(\frac{\alpha C^\gamma}{\gamma A^\alpha} x^{\alpha-1}\right)^{\gamma/(\gamma-1)} = 1 \qquad (22)$$

Because the left part of (22) is an increasing function of $x$ this equation has the unique solution and thus the rectangular parallelepiped of the greatest sum of the edges is unique, providing that it exists□

In the special case when $\alpha = \beta = \gamma$ we have a better result.

**Proposition 10.** There exists the unique parallelepiped of the greatest sum of the edges inscribed in the region bounded by the surface $\dfrac{|x|^\alpha}{A^\alpha}+\dfrac{|y|^\alpha}{B^\alpha}+\dfrac{|z|^\alpha}{C^\alpha}=1, \ \alpha>1.$

Proof. Let $x_0$ be the unique solution of equation (22) and let $y_0, z_0$ be the corresponding values obtained from (22). We claim that the rectangular parallelepiped with the vertex $(x_0, y_0, z_0)$ in the first octant has the greatest sum of the edges. From equations (21) and the equation of the surface follows that $x_0+y_0+z_0=\alpha\lambda$. We need to prove that the greatest value of the sum $x+y+z$ on the boundary of the surface is strictly smaller than $\alpha\lambda=x_0+y_0+z_0$. Consider first the point $(A,0,0)$. From the Lagrange equation $\dfrac{\alpha\lambda x_0^{\alpha-1}}{A^\alpha}=1$ we get $\alpha\lambda=A\dfrac{A^{\alpha-1}}{x_0^{\alpha-1}}>A$. Similarly, the global maximum of $x+y+z$ cannot be attained at $(0,B,0)$ or at $(0,0,C)$. Assume now that the maximum is attained at some point $(x_1, y_1, 0)$. At this point we have $\dfrac{\alpha\lambda_1 x_1^{\alpha-1}}{A^\alpha}=\dfrac{\alpha\lambda_1 y_1^{\alpha-1}}{B^\alpha}$ and $x_1+y_1=\alpha\lambda_1$. We have that either $x_1>x_0$ or $y_1>y_0$. In the first case $\alpha\lambda_1=\dfrac{A^\alpha}{x_1^{\alpha-1}}<\dfrac{A^\alpha}{x_0^{\alpha-1}}=\alpha\lambda$. The second case is similar.□

**Problem 8.** (1) Describe all the triples $(\alpha,\beta,\gamma)$ such that there exist rectangular parallelepipeds of the greatest and/or smallest sum of the edges inscribed in the surface (8).

(2) . Let S be a compact closed smooth surface in $\mathbf{R}^3$ symmetric about coordinate planes. Assume that the solid bounded by S is strictly convex (the intersection of every plane tangent to S with S is a singleton). Is it true that the rectangular parallelepiped (with faces parallel to coordinate planes) of the greatest sum of the edges inscribed in S exists and is unique?

## Part 2D. Inscribing a rectangular parallelepiped with the greatest ratios $\frac{V}{S^{3/2}}, \frac{V}{L^3}, \frac{S}{L^2}$ in a surface symmetric about coordinate planes.

In the title of this subsection $V$ is the volume of a rectangular parallelepiped, $S$ is its surface area, and $L$ is the sum of its edges.

**Proposition 11.** Assume that $f, g, h$ are functions continuous and strictly increasing on a closed interval $[0, A]$, continuously differentiable on $(0, A)$ and such that $f(0) = g(0) = h(0) = 0$. Let $S$ be the surface defined by the equation $f(|x|) + g(|y|) + h(|z|) = 1$. The maximum values of the ratios $\frac{V}{S^{3/2}}, \frac{V}{L^3}, \frac{S}{L^2}$ are attained in the case when the rectangular parallelepiped inscribed in $S$ is a cube.

*Proof.* We will prove our claim only for the ratio $\frac{V}{S^{3/2}}$, the other two cases can be considered similarly. Our problem is equivalent to the following one: minimize $F(x, y, z) = \frac{(xy + xz + yz)^3}{(xyz)^2}$ subject to $f(x) + g(y) + h(z) = 1$. The method of Lagrange multipliers provides $\frac{\partial F}{\partial x} = \lambda f'(x), \frac{\partial F}{\partial y} = \lambda g'(y), \frac{\partial F}{\partial z} = \lambda h(x)$. Therefore,

$$x\frac{\partial F}{\partial x} + y\frac{\partial F}{\partial y} + z\frac{\partial F}{\partial z} = \lambda(xf'(x) + yg'(y) + zh'(z)) \quad . \quad (23)$$

It is immediate to see that $F(x, y, z) = G\left(\frac{x}{y}, \frac{y}{z}\right)$ and therefore, (see e.g. [PZM])

$x\frac{\partial F}{\partial x} + y\frac{\partial F}{\partial y} + z\frac{\partial F}{\partial z} = 0$. Hence, $\lambda = 0$ and therefore, $\frac{\partial F}{\partial x} = \frac{\partial F}{\partial y} = \frac{\partial F}{\partial z} = 0$.

Straightforward computations show that this system of equations is equivalent to

$$xy + xz - 2yz = 0$$
$$xz + yz - 2xy = 0$$
$$xy + yz - 2xz = 0$$

Thus, $x = y = z = c$, where $c$ is defined in the unique way from the equation $\frac{c^\alpha}{A^\alpha} + \frac{c^\beta}{B^\beta} + \frac{c^\gamma}{C^\gamma} = 1$. □

To conclude this paper, we state our last problem.

**Problem 9.** What are possible variants of the results and problems discussed in this paper in the case of higher dimensions?

Arkady Kitover, Community College of Philadelphia, akitover@ccp.edu

Mehmet Orhon, University of New Hampshire, mehmetzorhon@gmail.com